\documentclass[runningheads,envcountsect,envcountsame]{llncs}

\usepackage[utf8]{inputenc}
\usepackage{amsmath, amsfonts, amssymb}
\usepackage{tikz}
\usetikzlibrary{automata, positioning, arrows}

\newcommand{\smat}[1]{{\small \arraycolsep=0.3\arraycolsep \renewcommand{\arraystretch}{0.8} \ensuremath{\begin{pmatrix}#1\end{pmatrix}} \renewcommand{\arraystretch}{1}}}
\title{Prefix palindromic length of the Sierpinski word}
\author{D. Bulgakova, A. Frid\inst{1}, J. Scanvic\inst{2}}
\institute{Aix Marseille Univ, CNRS, Centrale Marseille, I2M, Marseille, France,  \email{anna.frid@univ-amu.fr} \and Unité de Mathématiques Pures et Appliquées (UMR 5669),\\
École normale supérieure de Lyon / CNRS / Inria
}

\begin{document}

\maketitle

\begin{abstract}
The prefix palindromic length $p_{\mathbf{u}}(n)$ of an infinite word $\mathbf{u}$ is the minimal number of concatenated palindromes needed to express the prefix of length $n$ of $\mathbf{u}$. This function is surprisingly difficult to study; in particular, the conjecture that $p_{\mathbf{u}}(n)$ can be bounded only if $\mathbf{u}$ is ultimately periodic is open since 2013. A more recent conjecture concerns the prefix palindromic length of the period doubling word: it seems that it is not $2$-regular, and if it is true, this would give a rare if not unique example of a non-regular function of a $2$-automatic word.

For some other $k$-automatic words, however, the prefix palindromic length is known to be $k$-regular. Here we add to the list of those words the Sierpinski word $\mathbf{s}$ and give a complete description of $p_{\mathbf{s}}(n)$.

\end{abstract}

\section{Introduction}
A palindrome is a word which does not change when read from left to right and from right to left, like $rotator$ or $abbaaaabba$. In this paper, we continue to study decompositions of words over a finite alphabet to a minimal number of palindromes: for example, for the word $w=ababbaabbbaaa$, this number is equal to $4$, since we can factorize this word as $(aba)(bbaabb)(b)(aaa)$ or as $(aba)(bb)(aabbbaa)(a)$, but cannot manage with  less than four palindromes. So, we can write that the {\it palindromic length} of $w$, denoted as PL$(w)$, is equal to $4$.

In 2013, Puzynina, Zamboni and the second author \cite{fpz} conjectured that if the palindromic length of factors of an infinite word $\mathbf{u}$ is bounded, then the word
$\mathbf{u}$ is ultimately periodic. This conjecture remains open despite a partial solution in the initial paper \cite{fpz} and later particular results \cite{br,frid,ruk}.
Saarela \cite{saarela} proved that the conjecture is equivalent to the same statement about prefixes, not all factors, of $\mathbf{u}$. His result makes reasonable to consider the {\it prefix palindromic length} $p_\mathbf{u}(n)$, which is also denoted as PPL$_\mathbf{u}(n)$ in previous papers. This function of an infinite word $\mathbf{u}$ and of $n \geq 0$, equal to the palindromic length of the prefix of length $n$ of $\mathbf{u}$ is thus conjectured to be unbounded for every word which is not ultimately periodic.

A natural exercise on every new function of an infinite word is to compute or to estimate it for classical examples like the Thue-Morse word and the Fibonacci word. The first of these problems appears to be not too complicated: the prefix palindromic length of the Thue-Morse word, which is $2$-automatic, appears to be $2$-regular, and its first differences are described as a fixed point of a $4$-uniform morphism \cite{tm}. At the same time, the question on the Fibonacci word has not been solved, moreover, it seems that its prefix palindromic length is even not Fibonacci-regular \cite{flp}.

Since these first exercises, a progress has been made in computing the prefix palindromic length of some more known words, including the Rudin-Shapiro word, the paperfolding word \cite{flp} and the Zimin word \cite{liZ}. Moreover, it has been proved that for every $k$-automatic word containing a finite number of distinct palindromes, the prefix palindromic length is $k$-regular \cite{flp}. But the most intriguing are the results of computational experiments suggesting that for example, for the period-doubling word, which is the fixed point of the morphism $a \to ab, b \to aa$, the prefix palindromic length is {\it not} a $2$-regular sequence \cite{flp}. At the moment, this is the second challenging conjecture on the prefix palindromic length, since normally, all reasonable functions of $k$-automatic words are $k$-regular.

Unable to solve any of the big conjectures, we continue collecting examples when the prefix palindromic length is predictably regular.  Here we prove it for the Sierpinski word, the $3$-automatic fixed point of the morphism $\varphi: a \to aba,b \to bbb $. The fact that its prefix palindromic length is unbounded was proved already in the initial paper \cite{fpz}. The first morphic description of that function was conjectured in the Master thesis of Enzo Laborde \cite{laborde}, but here we find a simpler one, which yet requires several pages of proofs. We have also checked the results with the Walnut software \cite{walnut}.

A possible continuation of this research is to find a larger class of $k$-automatic words with $k$-regular prefix palindromic length. It could help to extract properties of automatic words which prevent the function to be regular.

The result can also be generalized to all morphisms of the form $a\to ab^{n-2}a, b\to b^n$ for $n\geq 3$, even though we do not include this result to this text.

\section{Definitions, notation, known results}
From now on, ${\bf s}=s[1]s[2]\cdots s[n]\cdots$ denotes the Sierpinski word, or the Cantor word
\[ababbbababbbbbbbbbababbbaba b^{27}\cdots,\]
defined as the fixed point starting with $a$ of the morphism $$\varphi: \begin{cases} a \to aba,\\ b \to bbb. \end{cases}$$ 
Here $s[i] \in \{a,b\}$ for all $i \geq 1$. Clearly, for every $k$, the Sierpinski word starts with the palindrome $\varphi^k(a)=\varphi^{k-1}(a)b^{3^{k-1}}\varphi^{k-1}(a)$. A factor $s[i]s[i+1]\cdots s[j]$ can also be denoted as ${\bf s}[i..j]$. The mirror image of a finite word $u$ is denoted by $\tilde{u}$; so, $u$ is a {\it palindrome} if $u=\tilde{u}$.

In what follows, PL$(u)$ denotes the palindromic length of a finite word $u$, that is, the minimal number of palindromes such that $u$ is their concatenation. The prefix palindromic length of ${\bf s}$ is denoted by $p_{\mathbf s}(n)$ or $p(n)$ for short: $p(n)=$PL$(s[1..n])$. A decomposition of a word $u$ to palindromes is {\it optimal} if the number of palindromes in it is minimal possible, that is, equal to PL$(u)$.

The {\it position} $n$ in a word is the position between its symbols numbered $n$ and $n+1$; it should not be confused with the symbol number $n$.

One of important general results on palindromic length is the following inequality, which we refer below as Saarela's inequality \cite[Lemma 6]{saarela}: for all words $u$, $v$ we have
\[|\mbox{PL}(u)-\mbox{PL}(v)|\leq \mbox{PL}(uv).\]
This result is especially useful when one of words $u$, $v$ or $uv$ is a palindrome and thus its palindromic length is equal to $1$. If $u$ is a prefix of a given infinite word ${\bf u}$ of length $n$, and $v$ is its next letter, it also immediately implies that 
\[|p_{\bf u}(n+1)-p_{\bf u}(n)|\leq 1,\]
meaning that the first differences of the prefix palindromic length of a word can be equal only to $-1$, $0$, or $1$.

As the name suggests, an infinite word ${\bf u}$ is called {\it $k$-automatic} if
there exists a deterministic finite automaton ${\bf A}$ such that every symbol $u[n]$
of ${\bf u}$ can be obtained as the output of ${\bf A}$ with the base-$k$ representation of $n$ as the input \cite{a_sh}. We will also need and use an equivalent definition of the same notion: a word ${\bf u}$ is $k$-automatic if and only if there exists a $k$-uniform morphism $\varphi: \Sigma^* \to \Sigma^*$ and a $1$-uniform morphism (or {\it coding}) $c: \Sigma^*\to \Delta^*$ such that ${\bf u}=c(\varphi^{\infty}(a))$ for a symbol $a \in \Sigma$.  So, for example, the Sierpinski word is $3$-automatic since its morphism $\varphi$ is $3$-uniform, and the coding $c$ can be chosen as the trivial one, sending $a$ to $a$ and $b$ to $b$.

A generalization of the notion of a $k$-automatic word to sequences on $\mathbb Z$ is the notion of $k$-regular sequence: formally speaking, a $\mathbb Z$-valued sequence is \emph{$k$-regular} if the $\mathbb Z$-module generated by its $k$-kernel is finitely generated. Discussions and equivalent definitions of $k$-regular sequences can be found in Chapter 16 of Allouche and Shallit's monograph \cite{a_sh}; what we really need in this paper is the following lemma proven in \cite{flp} for the case when $p$ is the prefix palindromic length of an infinite word but true for every sequence with bounded first differences due to exactly the same arguments.

\begin{lemma}
A $\mathbb Z$-valued sequence $r(n)$ with bounded first differences $d_r(n)=r(n+1)-r(n)$ is $k$-regular if and only if the sequence $d_r$ is $k$-automatic.
\end{lemma}

Since the main object we study in this paper is a $3$-automatic word, we need some addition notation concerning ternary representations.

Let $X \subset \{{\bf 0},{\bf 1}, {\bf 2}\}^*$ be the language of ternary expansions of non-negative integers without leading zeros. The fact that $x \in X$ is the ternary expansion of $n$ will be denoted as $[x]_3=n$ and $(n)_3=x$. By a convention, we put $(0)_3=\varepsilon$, so, the ternary representation of $0$ is the empty string. In means that every non-empty representation starts with ${\bf 1}$ or ${\bf 2}$, so, $X=\{\varepsilon\}\cup 
\{{\bf 1},{\bf 2}\}\{{\bf 0},{\bf 1}, {\bf 2}\}^*$. Note that we write symbols of ternary strings in boldface to distinguish concatenated strings from multiplied numbers. 

When we consider ternary expansions  with leading zeros, we mention that they are just strings over $\{{\bf 0},{\bf 1}, {\bf 2}\}^*$, not always from $X$.

For every function $f(n)$, we also use the notation $f(x)$, where $x$ is a ternary expansion of $n$.  
Also, let $x$ be a ternary expansion of $n$, where $2\cdot3^{k-1}\leq n \leq 3^k$; then we denote the ternary expansion of $3^k-n$ without leading zeros by $\overline{x}$. 

For all $k$, we clearly have $\overline{\bf 10}^k=\varepsilon$. For $x={\bf 1}x'$, where $x'\notin {\bf 0}^*$, the function $\overline{x}$ is not defined, for any other $x\in X$ we have $\overline{{\bf 2}x}=\overline{x}$.

\section{Auxiliary functions $q_j(n)$}
To study the prefix palindromic length $p(n)$ of the Sierpinski word, we first define for every $j \geq 0$ an auxiliary function
\[q_j(n)=\mbox{PL}(b^j s[1..n]).\]
Clearly, $p(n)=q_0(n)$, but for what follows, we need to study these functions for all $j$.
\begin{proposition}\label{p:qj}
 The functions $q_j$ can be found as follows:
 \begin{itemize}
  \item $q_0(0)=0$; for $j>0$, we have $q_j(0)=1$;
  \item $q_0(1)=1$; for $j>0$, we have $q_j(1)=2$;
  \item for $3^k\leq n \leq 2\cdot3^k$, we have $q_j(n)=1$ if $n=3^k+j$ and $q_j(n)=2$ otherwise;
  \item for $2\cdot3^k\leq n \leq 3^{k+1}$ and $j\leq 3^k$, we have
  \[q_j(n)=1+\min(q_{3^k-j}(n-2\cdot3^k),q_j(3^{k+1}-n)),\]
  \item at last, for $2\cdot3^k\leq n \leq 3^{k+1}$ and $j> 3^k$, we have
  \[q_j(n)=\min_{m\leq 3^k} q_m(n)+1.\]
 \end{itemize}

\end{proposition}
{\sc Proof.} The first three cases are obvious. To observe the last case, it is sufficient to see that $s[1..n]$ does not contain $b^j$, so, every decomposition of $b^js[1..n]$ starts with a palindrome $b^{j-m}$ for some $m$ and continues by a decomposition of $b^ms[1..n]$; we may choose $m$ to be the most convenient.

It remains to consider $2\cdot3^k\leq n \leq 3^{k+1}$ and $j\leq 3^k$. Here we know that ${\bf s}[1..n]=\varphi^k(a) b^{3^k}w$ for a prefix $w$ of $\varphi^k(a)$; so, $w=s[1..n-2\cdot3^k]$.

{\bf Type 1.} Let us consider the best of decompositions $P_1\cdots P_l=b^js[1..n]$ such that there is a boundary between palindromes at a position contained in $s[3^{k}..2\cdot3^{k}+1]=ab^{3^{k}}a$. If these boundaries are several, consider the first of them, denoted $m$, so that the prefix cut here is $u=b^j\varphi^k(a)b^m$. If $m>j$, then PL$(u)=2$, and among the decompositions of $u$ to two palindromes, we may choose $u=(b^j \varphi^k(a)b^j)(b^{m-j})$, so, first position $m$ will be replaced by $j$. Now suppose that $m<j$; then we also have PL$(u)=2$. Suppose that every decomposition of $b^js[1..n]$ starting with the palindrome $b^j \varphi^k(a)b^j$ is worse and contains $r>l$ palindromes. It means that PL$(b^{3^k-j}w)=r-1\geq l$, whereas PL$(b^{3^k-m}w)=l-2$. But this is impossible by the Saarela's inequality, since the first of these words is the suffix of the second preceeded by one palindrome $b^{j-m}$.

So, anyway, we may take $j=m$, and see that $l=1+$PL$(b^{3^k-j}w)=1+q_{3^k-j}(n-2\cdot3^k)$.

{\bf Type 2.} Now consider the shortest among palindromic decompositions $Q_1\cdots Q_{l'}=b^js[1..n]$ such that the word $s[3^{k}..2\cdot3^{k}+1]=ab^{3^{k}}a$ is contained in one palindrome. Since inside there is the largest power of $b$ in the considered word, this palindrome is the central part of some $Q_i=s[m_1+1..m_2]$, where $m_1=3^{k+1}-m_2$. Let us show that among such shortest decompositions, we can choose one with $m_2=n$.

Suppose we cannot do it. Choose another $m_2<n$ and denote $u=s[m_2+1..n]$. Since the considered decomposition is the best of its type,  we have
\[l'=q_j(m_1)+1+\mbox{PL}(u)=\mbox{PL}(b^jv\tilde{u})+1+\mbox{PL}(u),\]
where $s[1..m_1]=v\tilde{u}$ (it ends with $\tilde{u}$ since $s[1..3^{k+1}]$ is a palindrome, and $v$ is a new notation for the remaining prefix).

We have conjectured that any decomposition with the last palindrome 
$\tilde{u}b^{3^{k}}u$ is not optimal, meaning that 
\[l'=\mbox{PL}(b^jv\tilde{u})+1+\mbox{PL}(u)<\mbox{PL}(b^j v)+1,\]
that is, 
\[\mbox{PL}(b^jv\tilde{u})<\mbox{PL}(b^j v)-\mbox{PL}(u).\]
But since $\mbox{PL}(u)=\mbox{PL}(\tilde{u})$, this contradicts to Saarela's inequality. So, we may choose 
$m_2=n$. Consequently, $i=l'$, $m_1=3^{k+1}-n$, and $Q_1\cdots Q_{l'-1}$ is the optimal decomposition of  $b^js[1..m_1]=b^js[1..3^{k+1}-n]$, so, $l'=1+q_j(3^{k+1}-n)$.

It remains to notice that the optimal decomposition is either of type 1, or of the complementing type 2: $q_j(n)=\min(l,l')= \min(1+q_{3^{k}-j}(n-2\cdot3^{k}), 1+q_j(3^{k+1}-n))$. \hfill $\Box$

\begin{proposition}\label{p:qjdiff}
 For every $k\geq 0$, $j \geq 0$ and every $n \leq 3^{k}$ we have
  \[|q_j(n)-q_j(3^{k}-n)|\leq 1.\]
\end{proposition}
{\sc Proof.} It is sufficient to see that $q_j(3^{k}-n)=$PL$(b^js[1..3^{k}-n])=$PL$(s[n+1..3^{k}]b^j)$, since the last two words are mirror images one of the other. Since $q_j(n)=$PL$(b^js[1..n])$ and $b^j s[1..n]s[n+1..3^k] b^j=b^js[1..3^k]b^j$ is a palindrome, the inequality is a particular case of Saarela's one. \hfill $\Box$.

\section{Function $q$ and its first differences}

In this section, we study another auxiliary function $\displaystyle q(n) = \min_j q_j(n)$.

\begin{proposition}\label{p:q}
 For every $n \in \mathbb N$ the following equalities hold:
 \begin{equation}\label{e_q} q(n)=\begin{cases}
         0, \mbox{~if~} n=0;\\
         1, \mbox{~if~} n=1 \mbox{~or~} 3^{k}\leq n \leq 2\cdot3^{k};\\
         \min(1+q(n-2\cdot3^{k}), 1+q(3^{k+1}-n)),  \mbox{~if~} 2\cdot3^{k}<n\leq 3^{k+1},
         
        \end{cases}\end{equation}
meaning also for $q$ as the function of $X$ that
\begin{equation}\label{e_q2}\begin{cases}
         q(\varepsilon)=0;\\
         q({\bf 1}y)=1 \mbox{~for all~} y\in \{{\bf 0},{\bf 1}, {\bf 2}\}^*;\\
          q({\bf 2}y)=1+\min(q(y), q(\overline{{\bf 2}y}))  \mbox{~for all~}y\in 
        \{{\bf 0},{\bf 1}, {\bf 2}\}^*.
         
        \end{cases}\end{equation}
\end{proposition}
\noindent {\sc Proof.} First of all, note that 
$q(2\cdot3^{k-1})=q({\bf 20}^{k-1})=1=1+q(0)$ for all $k>0$, so, both \eqref{e_q} and \eqref{e_q2} are true for such values. In all other cases, the two statements are equivalent, so, it is sufficient to prove \eqref{e_q}. In fact, it immediately follows from Proposition \ref{p:qj} when we take the minimum for all $j$. 
\hfill $\Box$

Here is a list of basic properties of the function $q$.

\begin{proposition}
 For every $k \geq 0$ and every $n \leq 3^k$, we have $|q(n)-q(3^{k}-n)|\leq 1$.
\end{proposition}
{\sc Proof.} Follows directly from the definition of $q(n)=\min_j q_j(n)$ and Proposition \ref{p:qjdiff}. Indeed, suppose that $q(n)\geq q(3^k-n)$ and $j$ is such that $q(3^k-n)=q_j(3^k-n)$. Clearly, $q(n)\leq q_j(n)$. So, $q(n)-q(3^k-n)\leq q_j(n)-q_j(3^k-n)\leq 1$. The case of $q(n)\leq q(3^k-n)$ is symmetric. \hfill $\Box$

\begin{corollary} For every $n$ such that $2\cdot3^{k}<n\leq 3^{k+1}$, we have
$|q(n-2\cdot3^{k})-q(3^{k+1}-n)|\leq 1$. 
\end{corollary}
{\sc Proof.} Follows immediately from the previous proposition and the fact that if $n'=n-2\cdot3^k$, then $3^{k+1}-n=3^{k}-n'$. \hfill $\Box$

The next several properties of $q(x)$, $x \in X$, follow from \eqref{e_q2} and are proved by the same type of induction.

\begin{lemma}\label{l:12}
 For every $x\in X \cap \{{\bf 0},{\bf 2}\}^*$, we have $q(x{\bf 1})=q(x{\bf 2})$. 
\end{lemma}
{\sc Proof.} We proceed by induction on the length of $x$. For $x=\varepsilon$, we have $q({\bf 1})=q({\bf 2})=1$, so the base of induction holds. Now consider $x={\bf 2}y$ where $y\in  \{{\bf 0},{\bf 2}\}^*$ (so that $y$ may contain leading zeros). We have 
$q(x{\bf 1})=q({\bf 2}y{\bf 1})=1+\min(q(y{\bf 1}), q(\overline{{\bf 2}y{\bf 1}}))$ and $q(x{\bf 2})=q({\bf 2}y{\bf 2})=1+\min(q(y{\bf 2}), q(\overline{{\bf 2}y{\bf 2}}))$. But $q(y{\bf 1})=q(y{\bf 2})$ by the induction hypothesis; moreover, by the same hypothesis, $q(\overline{{\bf 2}y{\bf 1}})=q(\overline{{\bf 2}y{\bf 2}})$ since 
 $\overline{{\bf 2}y{\bf 1}}=z{\bf 2}$ and $\overline{{\bf 2}y{\bf 2}}=z{\bf 1}$ for the same $z\in X$, where $z$ is shorter than $x$. \hfill $\Box$

\begin{lemma}\label{l:0}
 For all $x \in X$, we have $q(x{\bf 0})=q(x)$.
\end{lemma}
{\sc Proof.} If $x=\varepsilon$, there is nothing to prove. If $x={\bf 1} y$, then $q(x{\bf 0})=q(x)=1$. Now for $x={\bf 2}y$, we proceed by induction on the length of $x$. The base is given by previous cases and $x={\bf 2}$ giving $q({\bf 20})=q(6)=q(2)=q({\bf 2})=1$. For the induction step, consider $x={\bf 2}y$, where the statement is proven for $y$ (which may start with leading zeros). It is sufficient to combine the last case of \eqref{e_q2} with the induction hypothesis and the fact that $\overline{x{\bf 0}}=\overline{x}{\bf 0}$, so that $q(y)=q(y{\bf 0})$, $q(\overline{x})=q(\overline{x}{\bf 0})=q(\overline{x{\bf 0}})$. \hfill $\Box$

\begin{lemma}\label{l:11w}
 For every $x\in X \cap \{0,2\}^*$ and for every $w \in\{{\bf 0},{\bf 1},{\bf 2}\}^*$, we have $q(x{\bf 1}w)=q(x{\bf 1})$. 
\end{lemma}
{\sc Proof.} As above, we start from $x=\varepsilon$ giving $q({\bf 1}w)=q({\bf 1})=1$ and proceed by induction on the length of $x$: take $x={\bf 2}y$ and suppose that the lemma is true for all strings shorter than $x$. As above, it is sufficient to compare $q(y{\bf 1})$ with $q(y{\bf 1}w)$, which are equal by the induction hypothesis, and $q(\overline{x{\bf 1}})$ with $q(\overline{x{\bf 1}w})$. 
For the latter comparison, we have to consider two cases: if $w \in \{\bf 0\}^*$, then the equality holds due to the previous lemma. 
If $w$ contains a non-zero symbol, then denote $\overline{x{\bf 1}}$ as $t{\bf 2}$ (indeed, its last symbol is equal to ${\bf 2}$). Then $\overline{x{\bf 1}w}=t{\bf 1} w'$ for some $w'$; but we know by from Lemma \ref{l:12} that $q(t{\bf 2})=q(t{\bf 1})$ and from the induction hypothesis that $q(t{\bf 1})=q(t{\bf 1}w')$. So, $q(\overline{x{\bf 1}})=q(\overline{x{\bf 1}w})$ and thus $q(x{\bf 1}w)=q(x{\bf 1})$. \hfill $\Box$

Summarizing Lemmas \ref{l:12} and \ref{l:11w}, we observe the following
\begin{corollary}
For every $x \in X$ such that $x=y{\bf 1}z$, where $y \in \{{\bf 0},{\bf 2}\}^*$, we have $q(x)=q(y{\bf 2})$.
\end{corollary}

So, we can concentrate on ternary representations from $ \{{\bf 0}, {\bf 2} \}^*$. and, due to Lemma \ref{l:0} even on those of them that end with ${\bf 2}$. 

For such a representation, that is, for a finite word on the alphabet 
$ \{{\bf 0}, {\bf 2} \}$, let us call a {\it small group} a sequence of ${\bf 2}$s separated from other such sequences by one or several ${\bf 0}$s. In its turn, a {\it large group} is a word beginning and ending with ${\bf 2}$ that does not contain two consecutive ${\bf 0}$s and is separated from other such groups by at least two consecutive ${\bf 0}$s. A large group is {\it dense} if it contains two consecutive ${\bf 2}$s and {\it sparse} otherwise.

\begin{example}{\rm
The word ${\bf 22202000022000202002}$ contains six small groups and four large groups (${\bf 22202}, {\bf 22}, {\bf 202}, {\bf 2}$). The first two of these large groups are dense and the last two are sparse. }
\end{example}

\begin{theorem}\label{t:q}
For every $x \in X \cap \{{\bf 0}, {\bf 2}\}^* {\bf 2}$, 
\begin{enumerate}
\item $q(x)=q(\overline{x})$ if and only if the first large group of $x$ is sparse, that is, if and only if ${\bf 20}^{k-1}\leq [x]_3 \leq ({\bf 20})^{k/2}$;  otherwise $q(x)=q(\overline{x})+1$;
\item
 the value of $q(x)$ is equal to the number of small groups plus the number of dense large groups in $x$.
 \end{enumerate}
 
\end{theorem}

Continuing the example above, we see that $q({\bf 22202000022000202002})=6+2=8$. Moreover, $\overline{\bf{ 22202000022000202002}}={\bf 20222200222020221}$, due to Lemma \ref{l:12}, $q({\bf 20222200222020221})=q({\bf 20222200222020222})$, and 
the latter representation contains 5 small groups and two large groups, both of them dense,  so that $q({\bf 20222200222020222})=q({\bf 20222200222020221})=7$. It is predicted by the first part of the theorem since the first large group of 
the initial representation is dense.

\medskip
{\sc Proof of the theorem.} 
As above, we proceed by induction on the length of $x$, but this time we have to consider several cases and prove both parts of the theorem together.

For the base of induction, consider $x={\bf 20}^i$, $i \geq 0$. Clearly, $\overline{x}={\bf 10}^i$, $x$ contains one small group and no dense large groups, and $q(x)=q(\overline{x})=1$, so, both statements hold.

Now for the induction step consider $x={\bf 20}^i x'$, where $x'$ starts with ${\bf 2}$, and suppose that the theorem is proven for all strings shorter than $x$.

\medskip
{\bf Case $i=0$:} suppose that $x={\bf 22}y$, $y \in \{{\bf 0}, {\bf 2}\}^*$. Then $q(x)=1+\min(q({\bf 2}y),q(\overline{x}))$; but $\overline{x}=\overline{{\bf 2}y}$, and so $q(x)=1+\min(q({\bf 2}y),q(\overline{{\bf 2}y}))$. 

{\bf Subcase ``dense''.}
If ${\bf 2}y$ starts with a dense large group, then so does $x$, and the number of small and dense large groups in $x$ is the same as in ${\bf 2}y$. Also, by the induction hypothesis, $q({\bf 2}y)=q(\overline{{\bf 2}y})+1=q(\overline{x})+1$, and thus $q(x)=q({\bf 2}y)=q(\overline{x})+1$. Both statements hold for $x$.

{\bf Subcase ``sparse''.}
If ${\bf 2}y$ starts with a sparse large group, then it becomes dense in $x$; the number of small groups stays the same. So, we should prove that $q(x)=q(\overline{x})+1$  for  the first statement of the theorem and that $q(x)=q({\bf 2}y)+1$ for the second one. Indeed, by the induction hypothesis,  $q({\bf 2}y)=q(\overline{{\bf 2}y})=q(\overline{x})$, so \eqref{e_q2} gives no choice for $q(x)$.

\medskip
{\bf Case $i=1$:} suppose that $x={\bf 202}y$, $y \in \{{\bf 0}, {\bf 2}\}^*$. For the second part of the statement, we should prove that $q(x)=q({\bf 2}y)+1$, since the number of small groups has increased and the number of dense large groups has not. For the first part of the statement, we should prove that $q(x)$ and $q(\overline{x})$ are in the same relation as $q({\bf 2}y)$ and $q(\overline{{\bf 2}y})$. In any case, $q(x)=1+\min(q({\bf 2}y), q(\overline{x}))$.

\smallskip
{\bf Subcase ``dense''.} Suppose that ${\bf 2}y$ starts with a dense large group, that is, ${\bf 2}y=({\bf 20})^j{\bf 22}z$ for some $j \geq 0$ and $z \in \{{\bf 0},{\bf 2}\}^*$. 

{\bf Subsubcase ``$z\in {\bf 0}^*$''.} Suppose first that $z= {\bf 0}^l$ for some $l$; then $\overline{{\bf 2}y}=({\bf 20})^j{\bf 10}^l$ and $\overline{x}=({\bf 20})^{j+1}{\bf 10}^l$. By the induction hypothesis, $q(\overline{x})=j+2$, $q({\bf 2}y)=j+2$ and $q( \overline{{\bf 2}y})=j+1$; then, due to \eqref{e_q2}, $q(x)=j+3$ and both statements hold.

{\bf Subsubcase ``$z\notin {\bf 0}^*$, $j=0$''.} In this case, $x={\bf 202}y={\bf 202}^{m}\cdots$, $m \geq 2$, and $\overline{x}={\bf 20}^m\overline{{\bf 2}y}$. By the induction hypothesis, $q(\overline{x})=q(\overline{{\bf 2}y})+1$, since $\overline{x}$ contains just one more small group in the beginning; also, ${\bf 2}y$ starts with a dense large group and thus $q({\bf 2}y)=q(\overline{{\bf 2}y})+1$. So, $q(\overline{x})=q({\bf 2}y)$ and $q(x)$ has to be equal to any of them plus one. Both statements hold.

{\bf Subsubcase ``$z\notin {\bf 0}^*$, $j>0$''.} The proof repeats the previous case but with  $\overline{x}={\bf 20}\overline{{\bf 2}y}$. 

\smallskip
{\bf Subcase ``sparse''.} Suppose that ${\bf 2}y$ (and $x$) start with a sparse large group, that is, ${\bf 2}y={\bf 2}$, or ${\bf 2}y={\bf 20}$ (these two cases are easy to consider separately), or ${\bf 2}y=({\bf 20})^j0z$ and thus $x=({\bf 20})^{j+1}0z$ for some $j\geq 1$ and $z \in \{{\bf 0}, {\bf 2}\}^*$. As above, we have to consider separately the cases when $z$ belongs or not to ${\bf 0}^*$, but in both cases, we have $\overline{x}={\bf 20}\overline{{\bf 2}y}$. So, $q(\overline{x})=1+q(\overline{{\bf 2}y})$ by the induction hypothesis, since just one small group is added. Also by the induction hypothesis, $q({\bf 2}y)=q (\overline{{\bf 2}y})$, since ${\bf 2}y$ starts with a sparse large group. Combining these equalities with \eqref{e_q2}, we see that $q(x)=q({\bf 2}y)+1=q(\overline{x})$, which was to be proved.

\medskip
{\bf Case $i\geq 2$.} Here $x={\bf 20}^i{\bf 2}y$ for some $y \in \{{\bf 0}, {\bf 2}\}^*$, with $i \geq 2$. The large first group of $x$ is sparse (and small), so, for this case, we should prove that $q(x)=q(\overline{x})$ and $q(x)=q({\bf 2}y)+1$; the second fact follows from the first and \eqref{e_q2}.

\smallskip
{\bf Subcase ``sparse''.} Suppose that ${\bf 2}y$ starts with a sparse large group and in particular, $y$ is either empty or starts with ${\bf 0}$.

{\bf Subsubcase $y={\bf 0}^j$, $j\geq 0$.} In this case, $x={\bf 20}^i{\bf 20}^j$, so, $\overline{x}={\bf 2}^{i}{\bf 10}^j$, and $q(\overline{x})=2$ by the induction hypothesis; at the same time, $q({\bf 20}^j)=1$, so, \eqref{e_q2} gives $q(x)=2=q(\overline{x})$, which was to be proved.

{\bf Subsubcase $y\notin {\bf 0}^*$.} In this case, since the starting group is sparse, $y$ starts with ${\bf 0}$, and thus $\overline{x}= {\bf 2}^i {\bf 0} \overline{{\bf 2} y}$. In particular, $q(\overline{x})>q( \overline{{\bf 2} y})$, since $\overline{x}$ contains at least one more small group ${\bf 2}^i$ at the beginning. At the same time, by the induction hypothesis, $q( \overline{{\bf 2} y})=q({\bf 2}y)$, so, $q(x)=1+\min(q(\overline{x}),q({\bf 2}y))=q(\overline{x})$, which was to be proved.

\smallskip
{\bf Subcase ``dense''.} Suppose that ${\bf 2}y$ starts with a dense large group.

{\bf Subsubcase $y\in {\bf 2}^+{\bf 0}^*$.} In this case, $x={\bf 20}^i{\bf 2}^j {\bf 0}^l$ with some $j\geq 2$. So, $\overline{x}={\bf 2}^i {\bf 0}^{j-1} {\bf 10}^l$ and thus, since $i \geq 2$, $q(\overline{x})=3$ by the induction hypothesis. Due to \eqref{e_q2}, $q(x)=1+\min(q(\overline{x}),q({\bf 2}^j {\bf 0}^l))$; since $q({\bf 2}^j {\bf 0}^l)=2$, the statement holds.

{\bf Subsubcase $y\notin {\bf 2}^+{\bf 0}^*$.} In this case, ${\bf 2}y$ contains at least two small groups; if the first of them is ${\bf 2}^j$, $j>0$, then $\overline{x}={\bf 2}^i {\bf 0}^j \overline{{\bf 2}y}$. If $j=1$, then ${\bf 2}^2$ is situated somewhere later in the first large group of ${\bf 2} y$, and so $\overline{{\bf 2}y}$ starts with a sparse group. In this case, $q(\overline{x})=q(
\overline{{\bf 2}y})+1$ since the prefix ${\bf 2}^i$ of $\overline{x}$ adds both a small group and a dense large group. The same is true if $j>1$, since in this case,  ${\bf 2}^i$ is itself a new dense large group in $\overline{x}$. 

At the same time, by the induction hypothesis, $q({\bf 2}y)=q(\overline{{\bf 2}y})+1$ since the first group is dense. So, $q(\overline{x})=q({\bf 2}y)+1$, and it remains to use \eqref{e_q2} to prove both statements in this last case. \hfill $\Box$

\bigskip

The first part of the theorem above will be used later for the results on the prefix palindromic length. As for the second part, it gives a formula for the function $q$ and in particular allows to find its first differences $d_q(n)=q(n+1)-q(n)$. The following corollary of the theorem is straightforward.

\begin{corollary}
 For every $n \geq 0$ with $(n)_3=x$ we have
 {\rm
 \[d_ q(n)=\begin{cases}
         0, \mbox{~if~} x \mbox{~contains~} {\bf 1}; \mbox{~otherwise}\\
        1, \mbox{~if~} x \mbox{~ends by~} {\bf 0} \mbox{~directly preceeded by~} {\bf 0} \mbox{~or a sparse large group};\\
         -1, \mbox{~if~} x \mbox{~ends by~} {\bf 2} \mbox{~which is a part of a dense large group};\\
         0, \mbox{~in all other cases}.
         \end{cases}\]
 }       

\end{corollary}
As it follows from this formula, the sequence $d_q(n)$ is automatic and here is the corresponding automaton.

Here and below, when considering first differences, we sometimes prefer to write $ \text{\tt -}$ instead of $-1$, $\text{\tt +}$ instead of $1$, and \text{\tt 0} in typewriter font.

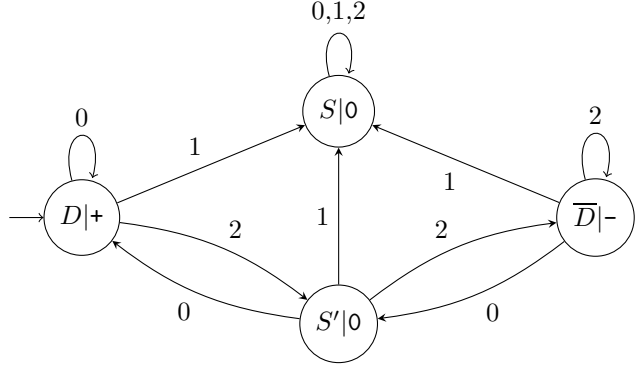
\begin{tikzpicture} [node distance = 2cm, on grid, auto]
 
\node (a) [state,initial,  initial text = {}] {$D|${\tt +}};
\node (inv1) [right=of a] {};
\node (s) [state, above right = of inv1] {$S|${\tt 0}};
\node (s1) [state, below right  = of inv1] {$S'|${\tt 0}};
\node (inv2) [below right = of s] {};
\node (a1) [state, right = of inv2] {$\overline{D}|${\tt -}};
 
\path [-stealth]
(a)  edge  [loop above] node {0} ()
      edge                       node {1} (s)
      edge [bend left=15]     node {2} (s1)
(s)  edge [loop above] node {0,1,2} ()
(s1) edge [bend left=15]   node {0} (a)
       edge                     node {1} (s)
       edge [bend left=15]   node {2} (a1)
(a1) edge [bend left=15]   node {0} (s1)
        edge                    node {1} (s)
        edge [loop above] node {2} ();
  
\end{tikzpicture} 

The choice of state names of this automaton will be clear from further constructions.

In its turn, this automaton is equivalent to the following morphic construction for the sequence $d_q$.
\begin{theorem} \label{t:dq} The sequence  $d_q$ is the 3-automatic word over the alphabet $\{ \text{\tt -}, \text{\tt 0}, \text{\tt +} \}$ given as follows:
 $$d_q=\gamma(\delta^{\infty}(D)),$$ where the morphism $\delta: \{D,S,S',\overline{D}\}^* \to \{D,S,S',\overline{D}\}^*$ is defined by
\[\begin{cases}
   \delta(D)=DSS',\\
   \delta(S)=SSS,\\
   \delta(S')=DS\overline{D},\\
   \delta(\overline{D})=S'S\overline{D},
  \end{cases}\]
  and the coding $\gamma: \{D,S,S',\overline{D}\}^* \to \{ \text{\tt -}, \text{\tt 0}, \text{\tt +} \}^*$ is given by $\gamma(D)=${\tt +}, $\gamma(S)=\gamma(S')=${\tt 0}, $\gamma(\overline{D})=${\tt -}.
\end{theorem}

\section{Difference between $p(n)$ and $q(n)$}
Now, after a study of the auxiliary function $q$, we return to the initial goal: the prefix palindromic length $p(n)$ of the Sierpinski word.

\begin{proposition} For every $n \geq 0$, the following holds.
 \[p(n)=\begin{cases}
         0, \mbox{~if~} n=0;\\
         1, \mbox{~if~} n=1;\\ 
         2, \mbox{~if~}  3^{k}<n\leq 2\cdot3^{k};\\
         \min(2+q(n-2\cdot3^{k}), 1+p(3^{k+1}-n)),  \mbox{~if~} 2\cdot3^{k}<n\leq 3^{k+1}.
         
        \end{cases}\]
Equivalently, this formula can be written as
\begin{equation}\label{e_p2}\begin{cases}
         p(\varepsilon)=0;\\
         p({\bf 10}^k)=1 \mbox{~for all~} k;\\
         p({\bf 1}y)=2 \mbox{~for all~} y \in \{{\bf 0}, {\bf 1}, {\bf 2} \}^*\backslash {\bf 0}^*;\\
          p({\bf 2}y)=1+\min(1+q(y), p(\overline{{\bf 2}y}))  \mbox{~for all~}y\in 
        \{{\bf 0},{\bf 1}, {\bf 2}\}^*.
         
        \end{cases}\end{equation}
\end{proposition}
{\sc Proof.} It is not difficult to see that the two statements are equivalent and that the first three lines hold. To prove the last equality note that $p(n)=q_0(n)$ where $q_0(n)=\min(1+q_{3^{k}}(n-2\cdot3^{k}),1-q_0(3^{k+1}-n))$ from Proposition \ref{p:qj}. From the last case of the same proposition, $q_{3^{k}}(n-2\cdot3^{k})=1+q(n-2\cdot3^k)$, so the equality follows. \hfill $\Box$

\begin{proposition}\label{p:3k-n}
 For every $n \geq 0$ such that $2\cdot3^k\leq n \leq 3^{k+1}$, the equality $p(n)=q(n)$ holds if and only if $p(3^{k+1}-n)=q(3^{k+1}-n)<q(n)$. Otherwise $p(n)=q(n)+1$.
\end{proposition}
{\sc Proof.} For the edge values, we easily check that $q(2\cdot3^k)=1<2=p(2\cdot3^k)$, and $q(3^{k+1}-2\cdot3^k)=q(2\cdot3^k)$, so that the condition does not hold; on the other hand, $q(3^{k+1})=p(3^{k+1})=1$, and the condition holds. For other values, from the previous results, we have
\[q(n)= \min(1+q(n-2\cdot3^{k}), 1+q(3^{k+1}-n)), \]
\[p(n)= \min(2+q(n-2\cdot3^{k}), 1+p(3^{k+1}-n)).\]
So, if $p(3^{k+1}-n)>q(3^{k+1}-n)$, then the values compared for $p(n)$ are just greater than the respective values compared for $q(n)$, and thus $p(n)>q(n)$. Moreover, suppose that $p(3^{k+1}-n)=q(3^{k+1}-n)$. If $q(3^{k+1}-n)=q(n)$, it immediately means that $q(n)=1+q(n-2\cdot3^{k})$ and $p(n)=1+q(3^{k+1}-n)=2+q(n-2\cdot3^k)>q(n)$. On the other hand, if $q(3^{k+1}-n)<q(n)$, then $q(n)=1+q(3^{k+1}-n)\leq 1+q(n-2\cdot3^k)$, so, $1+q(3^{k+1}-n)=1+p(3^{k+1}-n)< 2+q(n-2\cdot3^k)$ and thus $p(n)=1+q(3^{k+1}-n)=q(n)$. The equivalence is established.
\hfill $\Box$

The following statement is a direct corollary of the previous proposition and the first part of Theorem \ref{t:q}.

\begin{proposition}
 For every $x \in X$, we have $p(x)=q(x)$ if and only if $x \in {\bf 10}^*$ or $x$ starts with ${\bf 2}$, $p(\overline{x})=q(\overline{x})$ and $[x]_3> ({\bf 20})^{|x|/2}$.
\end{proposition}

Now the following statement can be proven by a straightforward induction.

\begin{proposition}\label{pr_n}
 Let $S \subset X$ be the set of ternary decompositions $x$ such that $p(x)=q(x)$. Then
 \[S=\{\varepsilon \} \cup \{{\bf 10}^*\} \cup \{({\bf 22}^+{\bf 00}^+)^*.{\bf 22}^+. \{{\bf 0^*} \cup {\bf 0}^+{\bf 10}^*\}.\]
\end{proposition}
In other words, $p(n)=q(n)$ if and only if $n=0$, $n=3^k$ for some $k$, or the ternary decomposition of $n$ consists of blocks of at least two ${\bf 2}$s and at least two ${\bf 0}$s, possibly followed by one ${\bf 0}$ or at least one ${\bf 0}$ before ${\bf 10}^l$ for some $l$.

\smallskip
\noindent {\sc Proof of the Proposition \ref{pr_n}.} Denote by $S_k$ the set of decompositions from $S$ corresponding to numbers not exceeding $3^k$ and by $D_k$ the difference $S_{k}\backslash S_{k-1}$. Clearly, Then $S_0=\{\varepsilon,{\bf 1}\}$, $D_1 =\{{\bf 10} \}$, $D_2 =\{{\bf 22, 100} \}$. Now, let us proceed by induction on $k$ starting with this base. Due to Proposition \ref{p:3k-n} for every $k$ we should look for elements of $D_{k+1}$ among numbers of the form $3^{k+1}-m$, $(m)_3\in S_k$. By the induction hypothesis, the elements of $D_k$ are  ${\bf 10}^k$ and some decompositions of length $k$ starting with ${\bf 22}$. They correspond to the numbers $m$ from $2\cdot 3^{k-1}+2 \cdot 3^{k-2}$ to $3^{k}$. So, if $(m)_3 \in D_k$, then $3^{k+1}-m\leq 3^{k+1}-2\cdot3^k-2\cdot3^{k-1}<2\cdot3^k$, and due to Proposition \ref{p:3k-n}, $3^{k+1}-m \notin S$. So, 
$$D_{k+1}=\{3^{k+1}-m|(m)_3 \in S_{k-1}\}.$$
It remains to check by a simple case study (whether $(m)_3$ contains ${\bf 1}$ or not) that subtracting from $3^{k+1}$ numbers whose ternary decompositions are in $S_{k-1}$ gives exactly numbers with decompositions from $S$, as described in the assertion, of length $k+1$, plus $3^{k+1}$.  \hfill $\Box$

Note also that the above expression for $D_{k+1}$ implies that
\[|D_{k+1}|=|S_{k+1}|-|S_k|=|S_{k-1}|,\]
and thus we can easily prove that every $|S_k|$ is a Fibonacci number: $|S_k|=F_{k+3}$ (if we start with $F_0=0, F_1=F_2=1$).

\medskip
The above proposition characterizes the function $t(n)=p(n)-q(n)$ which is equal to $0$ if $(n)_3 \in S$ and to $1$ otherwise. It also allows to find precisely its first differences $d_t(n)=t(n+1)-t(n)$:

\begin{corollary}
 The first differences of the function $t(n)$ are
 \[d_t(n)=\begin{cases}
           0, \mbox{~if~} (n)_3\in S \mbox{~does not contain~} {\bf 1} \mbox{~and ends with } {\bf 00} \mbox{~or~} {\bf 22};
        \\
1, \mbox{~if~} (n)_3\in S  \mbox{~contains~} {\bf 1} \mbox{~or ends with~} {\bf 220};\\
-1, \mbox{~if~}(n)_3 \in ({\bf 22}^+{\bf 00}^+)^*.\{{\bf 2} \cup {\bf 2}^+{\bf 12}^*\};\\
0, \mbox{~in all other cases}.

          \end{cases}
\]
Here the first case corresponds to $t(n)=t(n+1)=0$ and the last case to $t(n)=t(n+1)=1$.
\end{corollary}
The corresponding automaton for $d_t(n)$ is depicted below.

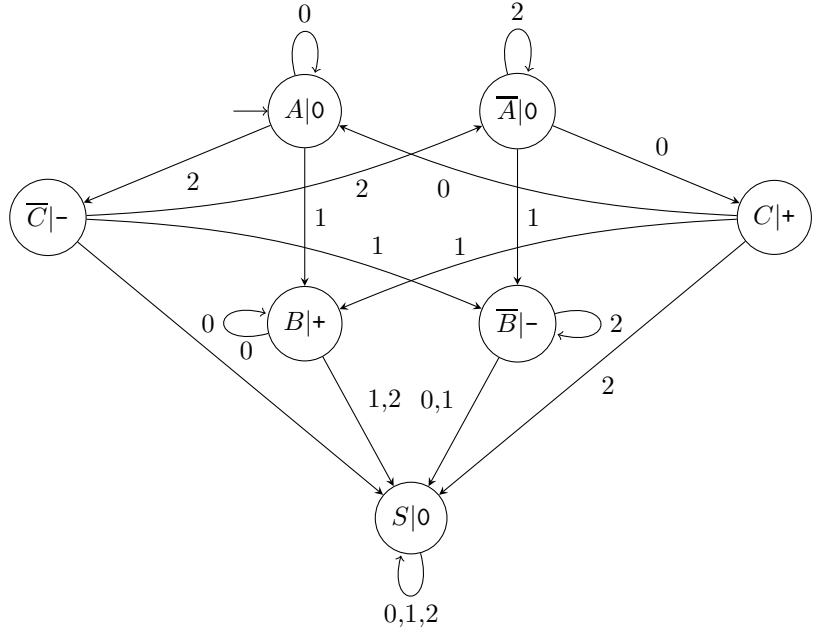
\begin{tikzpicture} [node distance = 2cm, on grid, auto]
\node (oC) [state] {$\overline{C}|${\tt -}};
\node (inv1) [right=of oC] {};
\node (A) [state,initial,  initial text = {},above right = of inv1] {$A|${\tt 0}};
\node (inv2) [below right = of A] {};
\node (oA) [state, above right = of inv2] {$\overline{A}|${\tt 0}};
\node (B) [state, below  left = of inv2] {$B|${\tt +}};
\node (oB) [state, below right = of inv2] {$\overline{B}|${\tt -}};
\node (inv3) [below right = of oA] {};
\node (C) [state, right = of inv3] {$C|${\tt +}};
\node (inv4) [below = of inv2] {};
\node (S) [state, below  = of inv4] {$S|${\tt 0}};

\path [-stealth]
(A)  edge  [loop above] node {0} ()
      edge                       node {1} (B)
      edge                       node {2} (oC)
(oA)  edge [loop above] node {2} ()
         edge                    node {1} (oB)
         edge                    node {0} (C)
(oC) edge   [bend right=10]                  node[below,pos=0.7] {2} (oA)
        edge    [bend left=10]                 node[pos=0.7] {1} (oB)
        edge                     node {0} (S)
(C) edge   [bend left=10]                  node[pos=0.7] {0} (A)
        edge    [bend right=10]                 node[above, pos=0.7] {1} (B)
        edge                     node {2} (S)
(B) edge [loop left] node{0} ()
      edge node{1,2} (S)
(oB) edge [loop right] node{2} ()
      edge node[above left]{0,1} (S)
(S) edge [loop below] node{0,1,2} ();
  
\end{tikzpicture} 

This automaton is equivalent to the following morphic construction for the sequence $d_t$.
\begin{theorem} \label{t:dt} The sequence  $d_t$ is the 3-automatic word over the alphabet $\{ \text{\tt -}, \text{\tt 0}, \text{\tt +} \}$ given as follows:
 $$d_t=\xi(\nu^{\infty}(A)),$$ where the morphism $\nu: \{A,B,C,\overline{A},\overline{B},\overline{C},S\}^* \to \{A,B,C,\overline{A},\overline{B},\overline{C},S\}^*$ is defined by
\[\begin{cases}
   \nu(A)=AB\overline{C},\\
   \nu(B)=BSS,\\
   \nu(C)=ABS,\\
   \nu(\overline{A})=C\overline{B}\overline{A},\\
   \nu(\overline{B})=SS\overline{B},\\
   \nu(\overline{C})=S\overline{B}\overline{A},\\
   \nu(S)=SSS,
  \end{cases}\]
  and the coding $\xi: \{A,B,C,\overline{A},\overline{B},\overline{C},S\}^* \to \{ \text{\tt -}, \text{\tt 0}, \text{\tt +} \}^*$ is given by $\xi(A)=\xi(\overline{A})=\xi(S)=${\tt 0}, $\xi(B)=\xi(C)=${\tt +}, $\xi(\overline{B})=\xi(\overline{C})=${\tt -}.
\end{theorem}

\section{First differences of $p(n)$}

By the definition of $t(n)$, the first differences of the function $p(n)$ are
\[d_p(n)=d_q(n)+d_t(n).\]
The functions $d_q(n)$ and $d_t(n)$ are completely described in Theorems \ref{t:dq} and \ref{t:dt} and by respective automata. It remains just to combine them, and one of the natural ways to do it is to define a new morphism $\psi =\smat{\delta\\ \nu}$ just as a direct product of $\delta$ and $\nu$ on the direct product of alphabets. We start with both starting symbols and get $\psi\smat{A\\D}=\smat{A\\D}\smat{B\\ S} \smat{\overline{C}\\S'}$; here the upper line is $\delta$ and the lower is $\nu$. Then we define $\psi$ on all the pairs of symbols that appeared, and continue this process while they continue to appear. We observe that only ten pairs appear in the fixed point of $\psi$ starting with $\smat{A\\D}$: the alphabet is  $\mathcal{A}=\left\{\smat{A\\D}, \smat{\overline{A}\\ \overline{D}}, \smat{B\\S}, \smat{\overline{B}\\S}, \smat{C\\S'}, \smat{\overline{C}\\S'},\smat{S\\D},\smat{S\\ \overline{D}}, \smat{S\\S}, \smat{S\\S'}\right\}$. Since we investigate the sum of the two first difference functions, each of these double letters is coded by $c\smat{X\\Y}=\gamma(X)+\xi(Y)$, where we recall that the symbols ${\tt -,0,+}$ are in fact numbers $-1,0,1$. So, for example, we have $c\smat{A\\D}=0+1=1$.

It remains to simplify the notation: the first six symbols of $\mathcal A$ can be denoted by just their upper letters, and the last four, starting with $S$, are defined by their lower letters. All this gives the following

\begin{theorem}\label{t:main}
 The sequence  $d_p$ of first differences of the prefix palindromic length of the Sierpinski word is the 3-automatic word over the alphabet $\{ \text{\tt -}, \text{\tt 0}, \text{\tt +} \}$ defined as 
 $$d_p=c(\psi^{\infty}(A)),$$ where the morphism $\psi: {\mathcal B}^* \to {\mathcal B}^*$, where ${\mathcal B}= \{A,B,C,D,\overline{A},\overline{B},\overline{C},\overline{D},S,S'\},$ is defined by
\[\begin{cases}
   \psi(A)=AB\overline{C},\\
   \psi(B)=BSS,\\
   \psi(C)=AB\overline{D},\\
   \psi(D)=DSS',\\
   \psi(\overline{A})=C\overline{B}\overline{A},\\
   \psi(\overline{B})=SS\overline{B},\\
   \psi(\overline{C})=D\overline{B}\overline{A},\\
   \psi(\overline{D})=S'S\overline{D},\\
   \psi(S)=SSS,\\
   \psi(S')=DS\overline{D},
  \end{cases}\]
  and the coding $c:{\mathcal B}^* \to \{ \text{\tt -}, \text{\tt 0}, \text{\tt +} \}^*$ is given by $c(A)=c(B)=c(C)=c(D)={\tt +}$, $c(\overline{A})=c(\overline{B})=c(\overline{C})=c(\overline{D})={\tt -}$, $c(S)=c(S')={\tt 0}$.
\end{theorem}
The corresponding DFAO is depicted below.

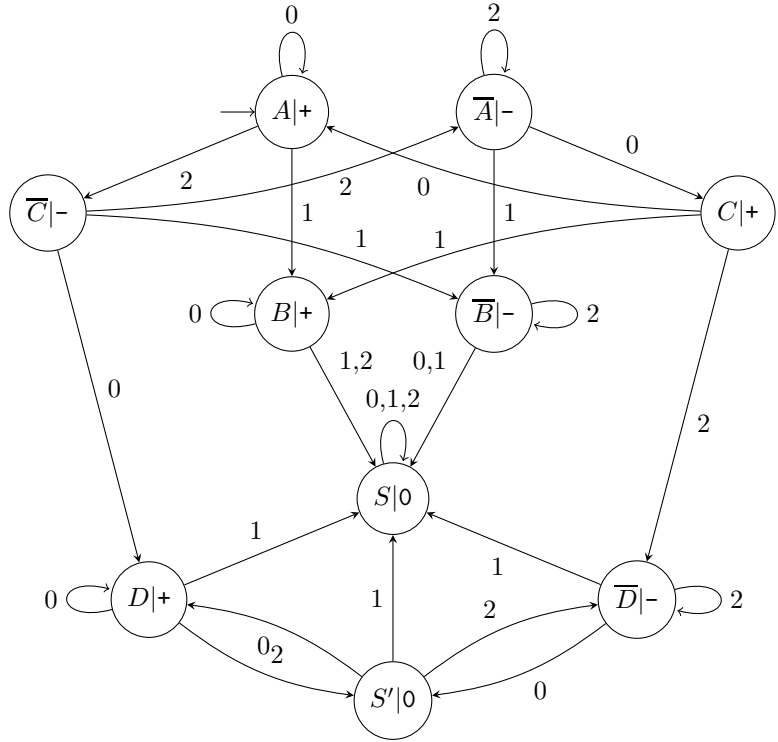
\begin{tikzpicture} [node distance = 1.9cm, on grid, auto]
\node (oC) [state] {$\overline{C}|${\tt -}};
\node (inv1) [right=of oC] {};
\node (A) [state,initial,  initial text = {},above right = of inv1] {$A|${\tt +}};
\node (inv2) [below right = of A] {};
\node (oA) [state, above right = of inv2] {$\overline{A}|${\tt -}};
\node (B) [state, below  left = of inv2] {$B|${\tt +}};
\node (oB) [state, below right = of inv2] {$\overline{B}|${\tt -}};
\node (inv3) [below right = of oA] {};
\node (C) [state, right = of inv3] {$C|${\tt +}};
\node (inv4) [below = of inv2] {};
\node (S) [state, below  = of inv4] {$S|${\tt 0}};
\node (inv5) [below left = of S] {};
\node (S1) [state, below right = of inv5] {$S'|${\tt 0}};
\node (D) [state, left = of inv5] {$D|${\tt +}};
\node (inv6) [below right = of S]{};
\node (oD) [state, right = of inv6] {$\overline{D}|${\tt -}};
 
\path [-stealth]
(A)  edge  [loop above] node {0} ()
      edge                       node {1} (B)
      edge                       node {2} (oC)
(oA)  edge [loop above] node {2} ()
         edge                    node {1} (oB)
         edge                    node {0} (C)
(oC) edge   [bend right=10]                  node[below,pos=0.7] {2} (oA)
        edge    [bend left=10]                 node[pos=0.7] {1} (oB)
        edge                     node {0} (D)
(C) edge   [bend left=10]                  node[pos=0.7] {0} (A)
        edge    [bend right=10]                 node[above, pos=0.7] {1} (B)
        edge                     node {2} (oD)
(B) edge [loop left] node{0} ()
      edge node[pos=0.3]{1,2} (S)
(oB) edge [loop right] node{2} ()
      edge node[above left, pos=0.3]{0,1} (S)
(S) edge [loop above] node{0,1,2} ()
(S1) edge [bend right=15]   node {0} (D)
       edge                     node {1} (S)
       edge [bend left=15]  node {2} (oD)
(oD) edge [bend left=15]   node {0} (S1)
        edge                    node {1} (S)
        edge [loop right] node {2} () 
(D) edge [loop left] node {0} ()
      edge                 node {1} (S)
      edge [bend right=15] node {2} (S1);
  
\end{tikzpicture} 


\vskip 0.2cm
We have proved that the first differences of the function $p_{\bf s}(n)$ are $3$-automatic and thus the function itself is $3$-regular.

\end{document}